\documentstyle{amsppt}
\TagsOnRight \headline={\tenrm\hss\folio\hss}
\magnification=\magstep1\nopagenumbers \parindent=24pt\vsize 22
true cm \hsize 16 true cm

\topmatter
\title On some residual properties of Baumslag -- Solitar  groups \endtitle
\author David Moldavanskii  \endauthor

\abstract A survey of results on the residual properties of
Baumslag -- Solitar groups which have been obtained to date.
\endabstract
\endtopmatter
\document

\centerline{\bf 1. Introduction}
\medskip

The Baumslag -- Solitar group (BS-group) is an one-relator group
with presentation
$$
G(m,n)=\langle a,\,b;\ a^{-1}b^{m}a=b^{n}\rangle,
$$
where $m$ and $n$ are non-zero integers. Note at once that since
groups $G(m,n)$, $G(n,m)$ and $G(-m,-n)$ are isomorphic we can
assume without loss of generality (and when it is convenient) that
integers $m$ and $n$ in the presentation of group $G(m,n)$
satisfy the condition $|n|\geqslant m>0$.

The family of groups $G(m,n)$ was introduced for consideration in
1962 in the paper of G.~Baumclag and D.~Solitar [2]. Just in this
family authors discovered the first examples of finitely generated
one-relator groups that are non-Hopfian (i.~e. are isomorphic to
some own proper quotient group) and therefore are not residual
finite; specifically, it was shown that the group $G(2,3)$ is
non-Hopfian. Thus, the supposition that every finitely generated
one-relator group is Hopfian turned out to be disproved. At that
time  some mathematicians believed that this assumption, as well
as the assumption of the residual finiteness of all one-relator
groups, is correct (perhaps because of the purely formal nearness
of one-relator groups and free groups). It should be noted also
that the properties of group $G(2,3)$ have given an answer to the
question of B.~H.~Neumann [17,~ p.~545] whether a 2-generator
non-Hopfian group can be defined by finite set of relations.

The study of properties of BS-groups became the permanent subject
of many investigations. This family of groups is of interest to
researchers, in particular, because some natural questions about
the properties of one-relator groups in the case of BS@-groups can
be answered in a more completed form than in the general case. For
example, the isomorphism problem for groups of this family is
trivial in view of following result (see [12]): groups $G(m,n)$
and $G(m',n')$, where $|n|\geqslant m>0$ and $|n'|\geqslant m'>0$,
are isomorphic if and only if $m=m'$ and $n=n'$. To a certain
extent the same is valid for problems about residual properties of
one-relator groups. This article is an extended version of [9] and
contains a survey of the results in this area that have been
received to date.

Some results are presented here with proofs. This generally
happens in cases where the relevant publication is inaccessible
now or (the new and more simple) proof has not been published.

Let us agree on the following terminology. If $\Cal K$ is a class
of groups then a group $G$ will be said to be {\it $\Cal
K$@-residual} if for any non-identity element $a\in G$ there
exists a homomorphism $\varphi$ of group $G$ onto some group from
class $\Cal K$ such that the image $a\varphi$ of $a$ is not equal
to identity. A group $G$ will be said to be {\it conjugacy $\Cal
K$@-separable} if for any  elements $a,b\in G$ that are not
conjugate in $G$ there exists a homomorphism $\varphi$ of group
$G$ onto some group $X$ from class $\Cal K$ such that the images
$a\varphi$ and $b\varphi$ of $a$ and $b$ are not conjugate in $X$.
Subgroup $H$ of group $G$ is said to be {\it $\Cal K$@-separable}
if for any element $g\in G\setminus H$ there exists a homomorphism
$\varphi$ of group $G$ onto some group from class $\Cal K$ such
that the image $g\varphi$ of element $g$ does not belong to image
$H\varphi$ of subgroup $H$. It is obvious that if a group is
conjugacy $\Cal K$@-separable then it is $\Cal K$@-residual and
group is $\Cal K$@-residual if and only if its identity
subgroup is $\Cal K$@-separable.

Let $\Cal F$ denote the class of all finite groups and if $p$ is a
prime number and $\pi$ is a set of prime numbers then let $\Cal
F_{p}$ and $\Cal F_{\pi}$ denote the class of all finite
$p$@-groups and the class of all finite $\pi$@-groups
respectively. It is clear that the property of $\Cal
F$@-residuality coincides with classical property of residuality
finite and the property of conjugacy $\Cal F$@-separability
coincides with classical property of conjugacy separability. Group
$G$ is said to be {\it subgroup separable} if all of its finitely
generated subgroups are $\Cal F$@-separable.\bigskip

\centerline{\bf 2. Residuality of BS-groups}
\medskip

The attempt to characterize $\Cal F$@-residual groups $G(m,n)$
made in  [2] was refined by S.~Meskin [8] as follows:

\proclaim{\indent Theorem 1} Group $G(m,n)$ is $\Cal F$@-residual
if and only if $($under the condition $|n|\geqslant m>0$$)$ either
$m=1$ or $|n|=m$.\endproclaim

The criterion of $\Cal F_p$@-residuality of groups $G(m,n)$ gives

\proclaim{\indent Theorem 2 {\rm (see [13, Theorem 3])}} For any
prime number $p$ group $G(m,n)$ $($where again it is supposed that
$|n|\geqslant m>0$$)$ is $\Cal F_{p}$@-residual if and only if
either $m=1$ and $n\equiv 1\pmod p$ or $|n|=m=p^{r}$ for some
$r\geqslant 0$ and also if $n=-m$ then $p=2$.\endproclaim

It makes sense to give a direct and quite elementary proofs of
these theorems. To do this we first note that any group $G(m,n)$
is an $HNN$@-extension with stable letter $a$ of infinite cyclic
base group $B$, generated by $b$, with associated subgroups
$B^{m}$ and $B^{n}$ that are generated by elements $b^{m}$ and
$b^{n}$ respectively. Secondly we introduce a family of finite
homomorphic images of group $G(1,n)$, namely, for arbitrary
positive integers $k$ and $l$ such that $n^{k}\equiv 1\pmod l$, we
set
$$
H_n(k,l)=\langle a,b;\ a^{-1}ba=b^{n},\,a^{k}=b^{l}=1\rangle.
$$
Since the order of automorphism of cyclic group $\langle b;\
b^{l}=1\rangle$ that is defined by the mapping $b\mapsto b^{n}$
divides the integer $k$, the group $H_n(k,l)$ is a split extension
of cyclic group $\langle b;\ b^{l}=1\rangle$ by cyclic group
$\langle a;\ a^{k}=1\rangle$. Hence, the order of group $H_n(k,l)$
is $kl$, orders of it's elements $a$ and $b$ are $k$ and $l$
respectively and any element $g\in H_n(k,l)$ can be uniquely
written in the form $g=a^{i}b^{j}$, where $0\leqslant i<k$ and
$0\leqslant j<l$.

Now, let $g$ be non-identity element of group $G(1,n)$. It is easy
to see (using relations $ba=ab^{n}$ ¨ $a^{-1}b=b^{n}a^{-1}$) that
element $g$ can be written as $g=a^{p}b^{s}a^{-q}$, where
$p,q\geqslant 0$, and therefore $g$ is conjugate to element
$a^{t}b^{s}$, where $t=p-q$. If $t\neq 0$ then the image of
element $g$ under obvious homomorphism of $G(1,n)$ onto infinite
cyclic group with generator $a$ is not equal to identity. If $t=0$
and hence $s\neq 0$ then the image of element $g$ in group
$H_n(k,l)$, where $l>0$ is chosen coprime to $n$ and not dividing
$s$ and $k=\varphi(l)$ is the value of the Euler function, is not
equal to identity.

Thus, the $\Cal F$@-residuality of any group $G(1,n)$ is proved.
Moreover, if for some prime number $p$ the congruence $n\equiv
1\pmod p$ is fulfilled then for any number $s>0$ we have
$n^{p^{s+1}}\equiv 1\pmod {p^{s}}$ and therefore the image of any
non-identity element $g\in G(1,n)$ in the suitable finite
$p$@-group $H_n(p^{s+1},p^{s})$ is not equal to identity.

If $|n|=m$, i.~e. $n=m\varepsilon$ for some $\varepsilon=\pm 1$,
then in group $G(m,m\varepsilon)$ subgroup $B^{m}$ is normal and
the quotient group $G(m,m\varepsilon)/B^{m}$ is the free product
of two cyclic groups, infinite and finite of order $m$. Therefore,
if non-identity element $g$ of group $G(m,m\varepsilon)$ does not
belong to subgroup $B^{m}$ then it's image in $\Cal F$@-residual
quotient group $G(m,m\varepsilon)/B^{m}$ is not equal to identity.
To consider the remaining case when $g=b^{ms}$ for some $s\neq 0$
let $\varphi$ be homomorphism of group $G(m,m\varepsilon)$ onto
group $G(1,\varepsilon)$ defined by identity mapping of
generators. Since  the group $G(1,\varepsilon)$ by above is $\Cal
F$@-residual and homomorphism $\varphi$ on subgroup $B$ acts
injectively the proof of $\Cal F$@-residuality of group
$G(m,m\varepsilon)$ is completed.

If $m=p^{r}$ for some prime number $p$ then the quotient group
$G(m,m\varepsilon)/B^{m}$ is $\Cal F_{p}$@-residual [4]. Moreover,
the group $G(1,1)$ is free Abelian and therefore is $\Cal
F_{p}$@-residual for any prime $p$. The group $G(1,-1)$ is $\Cal
F_{2}$@-residual since it's elements $a^{2}$ and $b$ generate free
Abelian normal subgroup of index 2.

Thus, the sufficiency of conditions in Theorems 1 and 2 is proved.
Let us show that these conditions are necessary.

If $|n|>m>1$ then element $b$ does not belong to subgroup $B^{m}$.
Also, if $d=(m,n)$ is the greatest common divisor of integers $m$
and $n$ then element $b^{d}$ does not belong to subgroup $B^{n}$.
Therefore the commutator $\left[ab^{d}a^{-1},\,b\right]$ is not
equal to 1 since it's expression
$ab^{-d}a^{-1}b^{-1}ab^{d}a^{-1}b$ is reduced in $HNN$@-extension
$G(m,n)$. On the other hand, turns out to be that this commutator
goes into the identity under any homomorphism of group $G(m,n)$
onto finite group. This assertion can be obtained from the
following observation:

\proclaim{\indent Proposition 1} Let elements $x$ and $y$ of a
group have the same finite order and let $x^{n}=y^{m}$ for some
integers $n$ and $m$. Then $\left[x^{d}, y\right]=1$ where
$d=(m,n)$ is the greatest common divisor of $m$ and
$n$.\endproclaim

Really, let $r=|x|=|y|$. Since $x^{n}=y^{m}$ we must have
$(r,n)=(r,m)$ and hence $(r,n)$ divides $d$. Consequently, there
exists an integer $s$ such that $ns\equiv d\pmod r$. Then
$x^{d}=x^{ns}=y^{ms}$ and therefore $\left[x^{d}, y\right]=1$ as
required.

Returning to the element $\left[ab^{d}a^{-1},\,b\right]$ of group
$G(m,n)$ it is sufficient to remark that if $\varphi$ is a
homomorphism of group $G(m,n)$ onto finite group then elements
$x=(aba^{-1})\varphi$ and $y=b\varphi$ satisfy the assumptions of
the Proposition 1.

So, the proof of Theorem 1 is complete. Now, let us suppose that
group $G(1,n)$ is $\Cal F_{p}$@-residual for some prime $p$. Then
there exists a homomorphism $\varphi$ of group $G(1,n)$ onto
finite $p$@-group $X$ such that $y=b\varphi\neq 1$. Let also
$x=a\varphi$. Since in group $G(1,n)$ for any number $k>0$ the
equality $a^{-k}ba^{k}=b^{n^{k}}$ holds, we have $n^{p^{r}}\equiv
1\pmod{p^{s}}$ where  $p^{r}$ is the order of element $x$ and
$p^{s}$ is the order of element $y$. Since $s>0$ this implies the
congruence $n^{p^{r}}\equiv 1\pmod{p}$. But as by Fermat Theorem
$n^{p-1}\equiv 1\pmod{p}$ and numbers $p^{r}$ and $p-1$ are
coprime we obtain the required congruence $n\equiv 1\pmod{p}$.

Next let us show that if group $G(m,m\varepsilon)$ is $\Cal
F_{p}$@-residual then $m$ is a $p$@-number. Indeed, otherwise
there exists a prime $q\neq p$ dividing $m$, $m=m_{1}q$. Then
$m>1$ and $m>m_{1}$ and therefore the commutator
$\left[a^{-1}b^{m_{1}}a,b\right]$ is a non-identity element of
group $G(m,m\varepsilon)$. On the other hand let $\varphi$ be a
homomorphism of group $G(m,m\varepsilon)$ onto finite $p$@-group
$X$, $x=a\varphi$ and $y=b\varphi$. Let also $p^{s}$ be the order
of element $y$. Since numbers $q$ and $p^{s}$ are coprime there
exists an integer $k$ such that $qk\equiv 1\pmod{p^{s}}$. Then
$x^{-1}y^{m_{1}}x=(x^{-1}y^{m}x)^{k}=y^{m\varepsilon k}$ and hence
$\left[a^{-1}b^{m_{1}}a,b\right]\varphi=1$.

Finally, we note that for any integer $k\geqslant 0$ in group
$G(m,m\varepsilon)$ the equality
$a^{-k}b^{m}a^{k}=b^{m\varepsilon^{k}}$ holds. Hence if
$\varepsilon=-1$ and if modulo some finite index normal subgroup
$N$ of group $G(m,m\varepsilon)$ the order $k$ of element $a$ is
an odd number then $b^{2m}\in N$. Therefore if a group $G(m,-m)$
is $\Cal F_{p}$@-residual then $p=2$ and Theorem 2 is
proved.\medskip

Theorems 1 and 2 can be generalized in the following way. Let
$\Cal K$ be again a class of groups and let for any group $G$ the
symbol $\sigma_{\Cal K}(G)$ denote the intersection of all normal
subgroups $N$ of group $G$ such that quotient group $G/N$ belongs
to $\Cal K$. It is clear that a group $G$ is $\Cal K$@-residual if
and only if $\sigma_{\Cal K}(G)$ coincides with identity subgroup.
Moreover, $\sigma_{\Cal K}(G)$ is the smallest normal subgroup of
$G$ the quotient group by which is $\Cal K$@-residual. If $\Cal
K=\Cal F$ or if $\Cal K=\Cal F_{p}$ then in place of $\sigma_{\Cal
K}(G)$ we shall write $\sigma(G)$ or $\sigma_p(G)$ respectively.

\proclaim{\indent Theorem 3 {\rm (see [10, Theorem 1])}} Let
$d=(m,n)$ be the greatest common divisor of integers $m$ and $n$.
Subgroup $\sigma\left(G(m,n)\right)$ coincides with the normal
closure in group $G(m,n)$ of the set of all commutators of form
$\left[a^{k}b^{d}a^{-k},\, b\right]$ where $k\in\Bbb Z$.
\endproclaim

\proclaim{\indent Theorem 4 {\rm (see [11])}} Let  $p$ be a prime
number and let $m=p^{r}m_{1}$ and $n=p^{s}n_{1}$ where
$r,s\geqslant 0$ and integers $m_{1}$ and $n_{1}$ are not divided
by $p$. Let also $d$ be the greatest common divisor of integers
$m_{1}$ ¨ $n_{1}$ and $m_{1}=du$ and $n_{1}=dv$. Then \roster
\item
if $r\neq s$ or if integers $m_{1}$ and $n_{1}$ are not congruent
modulo $p$ then subgroup $\sigma_{p}(G(m,n))$ coincides with the
normal closure in group $G(m,n)$ of element $b^{p^{t}}$ where
$t=\text{\rm min}\{r,s\}$;
\item
if $r=s$ and $m_{1}\equiv n_{1}\!\pmod p$ then subgroup
$\sigma_{p}(G(m,n))$ coincides with the normal closure in group
$G(m,n)$ of set consisting of element
$a^{-1}b^{p^{r}u}ab^{-p^{r}v}$ and of all commutators of form
$\bigl[a^{k}b^{p^{r}}a^{-k},b\bigr]\ (k\in \Bbb Z)$.
\endroster\endproclaim

It should be emphasize that in proofs of Theorems 3 and 4
criterions of $\Cal F$@-residu\-al\-ity and $\Cal
F_p$@-residuality of group $G(m,n)$ stated in Theorems 1 and 2 are
not used. Vice versa, Theorems 1 and 2 can be deduced from
Theorems 3 and 4 respectively.

To demonstrate this let me show, at first, how the sufficiency of
conditions in Theorem~1 for group $G(m,n)$ (where $|n|\geqslant
m>0$) to be $\Cal F$@-residual can be derived from Theorem 3. It
is well known (and easily to see) that if $m=1$ then the normal
closure in group $G(m,n)$ of element $b$ is the locally cyclic and
therefore Abelian group. Hence, all commutators of form
$\left[a^{k}b^{d}a^{-k},\, b\right]$ are equal to 1. If $|n|=m$
then  $d$, the greatest common divisor of integers $m$ and $n$, is
equal to $m$ and the defining relation of group $G(m,n)$ is of
form $a^{-1}b^{d}a=b^{d\varepsilon}$ for some $\varepsilon=\pm 1$.
Consequently, for any integer $k$ in group $G(m,n)$ we have the
equality $a^{k}b^{d}a^{-k}=b^{d\varepsilon^{k}}$ which implies
that again $\left[a^{k}b^{d}a^{-k},\, b\right]=1$. Thus, we see
that if either $m=1$ or $|n|=m$ then by Theorem~3 subgroup
$\sigma\left(G(m,n)\right)$ of group $G(m,n)$ is equal to identity
and therefore the group $G(m,n)$ is $\Cal F$@-residual.

Conversely,  if $|n|>m>1$ then, as was shown above, the commutator
$\left[ab^{d}a^{-1},\,b\right]$ is not equal to 1. Consequently,
Theorem 3 implies that subgroup $\sigma\left(G(m,n)\right)$ is not
equal to identity and therefore the group $G(m,n)$ is not $\Cal
F$@-residual.
\smallskip

Now, let us deduce Theorem 2 from Theorem 4.

Suppose that group $G(m,n)$ is $\Cal F_p$@-residual, i.~e.
$\sigma_{p}(G(m,n))$ coincides with identity subgroup. Since for
any $t\geqslant 0$ element $b^{p^{t}}$ differs from identity and
therefore does not belong to subgroup $\sigma_{p}(G(m,n))$, the
structure of this subgroup should be described in item (2) of
Theorem 4. Consequently, we see that (in notations from the
statement of Theorem~4) $r=s$ and $m_{1}\equiv n_{1}\!\pmod p$.
So, if $m=1$ and therefore $r=s=0$, $m_{1}=1$ and $n=n_{1}$, then
we obtain $n\equiv 1\!\pmod p$.

Next, we claim that if $m>1$ then $m_{1}=1=|n_{1}|$. Indeed, since
$\sigma_{p}(G(m,n))=1$ then by item (2) in group $G(m,n)$ all
commutators of form $\bigl[a^{k}b^{p^{r}}a^{-k},b\bigr]$ must be
equal to identity. But if $m_{1}>1$ then element $b^{p^{r}}$ does
not belong to subgroup $B^{m}$. Also, since $|n|>1$ element $b$
does not belong to subgroup $B^{n}$. Hence the expression
$$
\bigl[a^{-1}b^{p^{r}}a,b\bigr]=a^{-1}b^{-p^{r}}ab^{-1}a^{-1}b^{p^{r}}ab
$$
of commutator $\bigl[a^{-1}b^{p^{r}}a,b\bigr]$ is reduced in
$HNN$@-extension $G(m,n)$ and therefore this commutator cannot be
equal to identity. Similarly, assumption that $|n_{1}|>1$ implies
impossibility of equation $\bigl[ab^{p^{r}}a^{-1},b\bigr]=1$.

Thus, we have $m=p^{r}$ and $n=p^{r}\varepsilon$ for some
$\varepsilon=\pm 1$. Finally, if $\varepsilon=-1$ then the
congruence $m_{1}\equiv n_{1}\!\pmod p$ implies that $p=2$.

Conversely, if $m=1$ and $n\equiv 1\!\!\pmod p$ then $r=0$, $s=0$,
$m_{1}=1$ and $n_{1}=n$. Hence the congruence $m_{1}\equiv
n_{1}\!\pmod p$ is fulfilled. Therefore, in this case subgroup
$\sigma_{p}(G(m,n))$ is the normal closure in group $G(m,n)$ of
set of elements stated in item (2) of Theorem 4. As under $m=1$
the normal closure in group $G(m,n)$ of element $b$ is Abelian
group, all commutators in this set are equal to identity. Since in
this case we also have $p^{r}u=m$ and $p^{r}v=n$, element
$a^{-1}b^{p^{r}u}ab^{-p^{r}v}$ is equal to identity too.
Consequently, subgroup $\sigma_{p}(G(m,n))$ coincides with
identity, i.~e. group $G(m,n)$ is $\Cal F_p$@-residual.

If either $m=n=p^{r}$ or $m=2^{r}$ and $n=-2^{r}$ then subgroup
$\sigma_{p}(G(m,n))$ is again the normal closure in group $G(m,n)$
of set of elements stated in item (2) of Theorem 4 and it is clear
that all these elements are equal to identity. Thus, in these
cases group $G(m,n)$ is $\Cal F_p$@-residual and  $\Cal
F_2$@-residual respectively.\smallskip

Another way to generalize Theorems 1 and 2 consists of study of
conditions for group $G(m,n)$ to be $\Cal F_{\pi}$@-residual for
some (non-empty) set of prime numbers $\pi$. In paper [6] was
proved the

\proclaim{\indent Theorem 5 {\rm (see [6, Theorem 1])}}  Let $\pi$
be a set of prime numbers. Group $G(1,n)$ is $\Cal
F_{\pi}$@-residual if and only if there exists a $\pi$@-number
$s>1$ coprime to $n$ and such that the order modulo $s$ of integer
$n$  is a $\pi$@-number too.\endproclaim

The criterion in Theorem 2 for group $G(1,n)$ to be $\Cal
F_{p}$@-residual is a special case of Theorem 5.  Indeed, if group
$G(1,n)$ is $\Cal F_{p}$@-residual then by Theorem 5 we have
$n^{p^{t}}\equiv 1 \pmod{p^{r}}$ for some numbers $t$ and $r>0$.
Then $n^{p^{t}}\equiv 1 \pmod{p}$ and since by Fermat Theorem
$n^{p-1}\equiv 1 \pmod{p}$ it follows that $n\equiv 1 \pmod{p}$.
Conversely, if $n\equiv 1 \pmod{p}$ then the order modulo $p$ of
integer $n$  is equal to 1 and therefore is a $p$@-number.
Consequently, group $G(1,n)$ is $\Cal F_{p}$@-residual by Theorem
5.

Theorem 2 implies certainly that group $G(1,n)$ is $\Cal
F_{\pi}$@-residual if the set $\pi$ contains at least one prime
divisor of integer $n-1$. On the other hand this Theorem can be
applied also to prove the existence of 2-element set $\pi$ that
contains no numbers from $\pi(n-1)$ and such that group $G(1,n)$
is $\Cal F_{\pi}$@-residual.

\proclaim{\indent Corollary {\rm (see [6, Theorems 2 and 3])}} Let
$\pi=\{p,q\}$ be a set consisting of two prime numbers $p$ and $q$
such that $p<q$ and both $p$ and $q$ do not divide the integer
$n-1$. Then group $G(1,n)$ is $\Cal F_{\pi}$@-residual if and only
if $(n,q)=1$, $p$ divides $q-1$ and the order modulo $q$ of
integer $n$  is a $p$@-number. Moreover, if $|n|>1$ then for any
prime number $p$ that does not belong to set $\pi(n-1)$ there
exists a prime number $q>p$ such that $q\notin\pi(n-1)$ and group
$G(1,n)$ is $\Cal F_{\pi}$@-residual where
$\pi=\{p,q\}$.\endproclaim

These results (Theorem 5 and Corollary) allows us to describe some
sets  $\pi$ of primes such that group $G(1,n)$ is $\Cal
F_{\pi}$@-residual and is not $\Cal F_{\pi_1}$@-residual for any
proper subset $\pi_1$ of $\pi$. For example, the group $G(1,2)$ is
not $\Cal F_{p}$@-residual for any prime $p$ and any prime $p$ is
contained in some 2-element set $\pi$ which is minimal such that
group $G(1,2)$ is $\Cal F_{\pi}$@-residual. In addition, since the
integer 2 is a primitive root modulo 29, the set $\pi=\bigl\{2, 7,
29\bigr\}$ is minimal with the property that group $G(1,2)$ is
$\Cal F_{\pi}$@-residual.\smallskip

When $|n|=m$, the criterion of $\Cal F_{\pi}$@-residuality of
group $G(m,n)$ can be expressed in more complete form:

\proclaim{\indent Theorem 6 {\rm (see [20, Theorem 2])}} Let $\pi$
be a set of prime numbers. Group $G(m,m)$ is $\Cal
F_{\pi}$@-residual if and only if  $m$ is a $\pi$-number, and
group $G(m,-m)$ is $\Cal F_{\pi}$@-resi\-dual if and only if  $m$
is a $\pi$-number and $\pi$ contains the integer $2$.\endproclaim

We conclude this section with recent results of D.~Azarov [1]
about virtually residuality of BS-groups. Recall that for any
class of groups  $\Cal K$ a group $G$ is said to be virtually
$\Cal K$@-residual if it contains a finite index subgroup which is
$\Cal K$@-residual. It is obvious that if the class $\Cal K$
consists only of finite groups, then any virtually $\Cal
K$@-residual group is $\Cal F$@-residual.

\proclaim{\indent Theorem 7 {\rm (see [1, Theorem 1])}} A group
$G(1,n)$ is virtually $\Cal F_{p}$@-residual if and only if the
prime $p$ does not divide $n$. If $|n|=m$ then for any prime $p$
group $G(m,n)$ is virtually $\Cal F_{p}$@-residual.\endproclaim

\proclaim{\indent Theorem 8 {\rm (see [1, Theorem 2])}} For any
set $\pi$ of prime numbers group $G(m,n)$ is virtually $\Cal
F_{\pi}$@-residual if and only if it is virtually $\Cal
F_{p}$@-residual for some $p\in\pi$.\endproclaim\bigskip

\centerline{\bf 3. Conjugacy separability of BS-groups}
\medskip

As it was noted above, any conjugacy $\Cal F$@-separable group is
$\Cal F$@-residual. For BS@-groups converse is also true:

\proclaim{\indent Theorem 9} If group $G(m,n)$ is $\Cal
F$@-residual then it is conjugacy $\Cal F$@-separable.\endproclaim

Conjugacy $\Cal F$@-separability of groups $G(1,n)$ was proved in
[14]. This assertion is contained also in more general result that
was obtained in [19] and affirms that any descending
$HNN$@-extension of finitely generated Abelian group  is a
conjugacy $\Cal F$@-separable group.

Conjugacy $\Cal F$@-separability of groups $G(m,n)$ when $|n|=m$
can be deduced from the result of work [21] or from generalization
of it which was obtained in [18]. It should be also noted that
since under $n=m$ the center of group $G(m,n)$ is non-trivial, the
statement on conjugacy $\Cal F$@-separability of group $G(m,n)$ in
this case follows as well from Armstrong's theorem which states
that any one-relator group with  non-trivial center is conjugacy
$\Cal F$@-separable (see e.~g. [3]).

However, we shall show now that in the case $|n|=m$ the statement
on the conjugacy $\Cal F$@-separability of group $G(m,n)$ can be
easily proved having applied  ideas of M.~I.~Kargapolov [7] and
result of J.~Dayer [3]. We reproduce also the original proof of
conjugacy $\Cal F$@-separa\-bi\-lity of group $G(1,n)$ given in
[14].\smallskip

\demo{\indent The proof of Theorem 9 in the case $m=1$}

Suppose that  the coprime integers $n\neq\pm 1$ and $k>0$ are
fixed. Then integers  $r$ and $s$ will be said to be
$(n,k)$@-equivalent, if there exists a number $x\geqslant 0$ such
that the  congruence $n^{x}r\equiv s\pmod k$ holds; it is obvious
that this relation is indeed an equivalence. It will allow us to
give the necessary and sufficient conditions for certain elements
of groups $G(1,n)$ and $H_n(r,s)$ (introduced above)  to be
conjugate. For any number $t>0$ we set $u_{t}=|n^{t}-1|$.

\proclaim{\indent Proposition 2} For any integer $n\neq\pm 1$  the
following assertions are true: \roster

\item every element of group $G(1,n)$ is conjugate to element of form $a^{t}b^{r}$ for
suitable integers $t$ and $r$ where if the number $r$ is not $0$,
then it is not divisible by $n$;

\item if $t>0$ then elements $a^{t}b^{r}$ and $a^{t}b^{s}$ are conjugate
in group $G(1,n)$ if and only if the integers $r$ and $s$ are
$(n,u_t)$@-equivalent.\endroster\endproclaim

The verity of the first part of (1) was noted above (in the proof
of Theorem 1). If $r\neq 0$ and $r=nr_1$ then element $a^{t}b^{r}$
is conjugate to element $a(a^{t}b^{r})a^{-1}=a^{t}b^{r_1}$ of the
same form with $|r_1|<|r|$. So, the truth of the second part (1)
is also proved.

To prove (2) we first assume that the elements  $a^{t}b^{r}$ are
$a^{t}b^{s}$ are conjugate in group $G(1,n)$, i.~e.
$g^{-1}(a^{t}b^{r})g=a^{t}b^{s}$ for some $g\in G(1,n)$. Let, as
above, $g=a^{p}b^{v}a^{-q}$, where $p,q\geqslant 0$. Then $
b^{-v}a^{-p}(a^{t}b^{r})a^{p}b^{v}=a^{-q}(a^{t}b^{s})a^{q}$, and
therefore
$$a^{t}\cdot(a^{-t}ba^{t})^{-v}\cdot(a^{-p}ba^{p})^r\cdot
b^{v}=a^{t}\cdot(a^{-q}ba^{q})^s.$$ Hence
$b^{n^pr-(n^t-1)v}=b^{n^qs}$ and since the order of element $b$ is
infinite we have the equality $n^pr-(n^t-1)v=n^qs$ from which the
congruence $n^pr\equiv n^qs\pmod{|n^t-1|}$ follows. Therefore the
integers $r$ and $s$ are $(n,u_t)$@-equivalent.

Conversely, if for some integer $x$ the congruence $n^{x}r\equiv
s\pmod {u_t}$ is valid then for suitable integer $y$ we have
$n^{x}r=s+(n^{t}-1)y$. Hence
$$
(a^{x}b^{y})^{-1}(a^{t}b^{r})(a^{x}b^{y})=
a^{t}\cdot(a^{-t}ba^{t})^{-y}\cdot(a^{-x}ba^{x})^{r}\cdot
b^{y}=a^{t}b^{n^{x}r-(n^{t}-1)y}=a^{t}b^{s},
$$
and Proposition 2 is proved.

\proclaim{\indent Proposition 3} The elements $b^{r}$ and $b^{s}$
of group $H_n(p,q)$ are conjugate if and only if the integers $r$
and $s$ are $(n,q)$@-equivalent.\endproclaim

Indeed, for any element $g=a^{i}b^{j}$ of group $H_n(p,q)$ the
equality $g^{-1}b^{r}g=b^{s}$ is equivalent to equality
$a^{-i}b^{r}a^{i}=b^{s}$ which, in turn, can be rewritten in the
form $b^{n^ir}=b^{s}$. Thus, the elements  $b^{r}$ ¨ $b^{s}$ are
conjugate if and only if for some integer $i\geqslant 0$ the
congruence $n^ir\equiv s\pmod q$ holds.
\medskip

A crucial role in the proof of the assertion of Theorem 9 in the
case $m=1$ plays the following statement from elementary number
theory.

\proclaim{\indent Proposition 4} Let $n$ be an integer $\neq\pm
1$. Then for any integers $r$ and $s$, where $r\neq s$ and both
$r$ and $s$ are not divisible by $n$, there exists a number $t>0$
such that the exponential congruence  $n^{x}r\equiv s \pmod{u_{t}}
$ has no solution.\endproclaim

The proof of Proposition 4 will be given below, after we use it to
complete the proof of conjugacy $\Cal F$@-separability of groups
$G(1,n)$.

It is obvious that the (free Abelian) group $G(1,1)$ is conjugacy
$\Cal F$@-separable. The conjugacy $\Cal F$@-separability of group
$G(1,-1)$  follows  from the  result of S.~M.~Armstrong mentioned
above, since the center of group $G(1,-1)$ is non-trivial. So, we
can assume that $n\neq\pm 1$.

Let $f$ and $g$ be the non-conjugate elements of group $G(1,n)$.
By the item (1) of Proposition 2 we  may suppose that
$f=a^{t_1}b^{r}$ and $g=a^{t_2}b^{s}$ for some integers $t_1$,
$t_2$, $r$ and $s$ such that if any of numbers $r$ and $s$ is not
equal to $0$, then it is not divisible by~$n$. If $t_1\neq t_2$
then the images of elements $f$ and $g$ under the evident
homomorphism of group $G(1,n)$ onto some finite cyclic group are
distinct and therefore are non-conjugate. Thus, it remains to
consider the case when $f=a^{t}b^{r}$ and $g=a^{t}b^{s}$. Here we
can assume also (replacing, if it is necessary, elements $f$ and
$g$ by $f^{-1}$ and $g^{-1}$) that $t\geqslant 0$.

If $t>0$ then by item (2) of Proposition 2 the integers $r$ ¨ $s$
are not $(n,u_t)$@-equiva\-lent. Therefore, by Proposition 3 the
images $b^{r}$ and  $b^{s}$ of elements $f$ and $g$ under natural
homomorphism of group $G(1,n)$ onto finite group $H_n(t,u_t)$ are
not conjugate in this group.

Finally, let $f=b^{r}$ and $g=b^{s}$. Since the group $G(1,n)$ is
$\Cal F$@-residual we can assume that both integers $r$ and $s$
are not equal to 0 and therefore are not divisible by $n$.  Then
by Proposition 4 there exists a number $t>0$ such that numbers $r$
and $s$ are not $(n,u_t)$@-equivalent. Consequently, the images of
elements $f$ and $g$ under homomorphism of group $G(1,n)$ onto
finite group $H_n(t,u_t)$ are not conjugate in this group. So, the
conjugacy $\Cal F$@-separability of groups $G(1,n)$ is
proved.\medskip

Now proceed to the proof of Proposition 4. It states that for any
integer $n\neq\pm 1$ and for any integers $r$ and $s$, $r\neq s$,
that are  not divisible by $n$ there exists a number $t>0$ such
that the exponential congruence $$n^{x}r\equiv s \pmod{u_{t}}
\tag1$$ has no solutions. To prove this, let us consider two cases
depending on the sign of $n$.\smallskip

{\it Case} 1, $n>0$. We shall show that in this case there exists
an integer $t_{0}>0$ such that for any $t\geqslant t_{0}$ the
congruence (1) does not have solution.

Assuming (without loss of generality) that the integer $r$ is
positive, we can write it in the number system with base $n$:
$$
r=c_{0}n^{k}+c_{1}n^{k-1}+\cdots +c_{k-1}n+c_{k},
$$
where $k\geqslant 0$, $0\leqslant c_{i}<n$ for any $i=0,1,\dots ,
k$ and $c_{0}\neq 0$. Remark that, since $r$ is not divisible by
$n$, we have also $c_{k}\neq 0$.

Next, let $l$ be a positive integer and $R=n^{l}r$. Then
$$
R=d_{0}n^{k+l}+d_{1}n^{k+l-1}+\cdots +d_{k+l-1}n+d_{k+l},
$$
where of course
$$
d_{i}=\cases c_{i},&\text{if\ $0\leqslant i\leqslant k$,}\\
0,&\text{if\ $k+1\leqslant i \leqslant k+l$.}\endcases
$$
Further, for every $i=0,1,\dots , k$ let the symbol $r_{i}$ denote the number that is obtained
from number $r$ by cyclic permutation of digits  beginning with $c_{i}$; thus, $r_{0}=r$ and
for $i>0$
$$
r_{i}=c_{i}n^{k}+c_{i+1}n^{k-1}+\cdots +c_{k}n^{i}+c_{0}n^{i-1}+\cdots +c_{i-1}.
$$
Similarly, for every $i=0,1,\dots , k+l$ let the number $R_{i}$ be
obtained by cyclic permutation of digits of number $R$ beginning
with $d_{i}$. Thus, $R_{0}=R$ and if $i>0$
$$
R_{i}= \sum_{j=0}^{k+l-i}d_{i+j}n^{k+l-j} +\sum_{j=0}^{i-1}d_jn^{i-1-j}.
$$

One can  easily show that under $t=k+l+1$ for any $i=0,1,\dots , k+l$ we have the congruence
$$
n^{i}R\equiv R_{i}\pmod{u_{t}}.\tag2
$$
Moreover, it is not difficult to see that
$$
R_{i}=\cases n^{l}r, & \text{ if\ $i=0$,} \\
n^{l}r_{i}+p_{i}(1-n^{l}),& \text{ if\ $1\leqslant i\leqslant k$,} \\
n^{i-k-1}r,& \text{ if\ $k+1\leqslant i\leqslant k+l$,}\endcases\tag3
$$
where for $1\leqslant i\leqslant k$ $p_{i}= c_{0}n^{i-1}+c_{1}n^{i-2}+\cdots +c_{i-1}$.

Congruences (2) obviously imply  that any integer of form $n^{i}R$, $i\geqslant 0$, is
congruent modulo $u_{t}$ (where, recall, $t=k+l+1$) to one of numbers $R_{0}$, $R_{1}$,\dots ,
$R_{k+l}$. From this and from (3) it follows that the same holds also for any number of form
$n^{i}r$. Indeed, if $i\geqslant l$  this is evident as $n^{i}r=n^{i-l}R$. In the case
$0\leqslant i \leqslant l-1$ we set $j=i+k+1$. Then $k+1\leqslant j \leqslant k+l$ and
therefore by (3) we have $n^{i}r=n^{j-k-l}r=R_j$. Remark also that $0<R_{i}<n^{k+l+1}$ for any
$i=0,1,\dots , k+l$.

Now, if in the case when $s>0$ we choose the number $l$ such that $n^{l}>s$ then all numbers
$s$ and $R_{0}$, $R_{1}$,\dots , $R_{k+l}$ will  belong to complete system of (the smallest
non-negative) residues modulo $u_{t}$. In addition, number $s$ is not equal to any number
$R_{i}$ ($0\leqslant i \leqslant k+l$). Really, if $i=0$ or $k+1\leqslant i \leqslant k+l$ this
follows directly from (3) since $s$ is different from $r$ and is not divisible by $n$. If
$1\leqslant i\leqslant k$ then again by (3) we have
$$
R_{i}= n^{l}(r_{i}-p_{i}) +p_{i}>n^{l}(c_{i}n^{k}+c_{i+1}n^{k-1}+\cdots +c_{k}n^{i})\geqslant
n^{l+i}c_{k}\geqslant n^{l+i}>s.
$$

Thus, if $s>0$ and if we set $t_{0}=k+l_{0}+1$, where $n^{l_0}>s$, then for any $t>t_{0}$ the
congruence (1) does not have solution.

In the case when $s<0$ it is sufficient to show that there exists a number  $l_{0}>0$ such that
$$
R_{i}<(n^{k+l+1}-1)+s \qquad (0\leqslant i \leqslant k+l)
$$
for any $l\geqslant l_{0}$. Indeed, then all numbers $s$ and $R_{0}$, $R_{1}$,\dots , $R_{k+l}$
will  belong to complete system $\bigl\{y\bigm |s\leqslant y<u_{t}+s\bigr\}$ of residues modulo
$u_{t}$ with $s<0<R_{i}$.

It follows from (3) that
$$
n^{k+l+1}-R_{i}=\cases n^{l}(n^{k+1}-r), & \text{ if $i=0$,} \\
n^{l}(n^{k+1}-r_{i}+p_{i})-p_{i},& \text{ if $1\leqslant i\leqslant k$,} \endcases
$$
and if $k+1\leqslant i\leqslant k+l$, then $ n^{k+l+1}-R_{i}\geqslant n^{l-1}(n^{k+2}-r)$.
Since all numbers $n^{k+1}-r$, $n^{k+1}-r_{i}+p_{i}$, $n^{k+2}-r$ are positive the existence of
the required number $l_{0}$ is evident.\smallskip

{\it Case} 2, $n<0$. If the integers $r^2$ and $s^2$ are distinct then, since they are not
divisible by $n^2$, it follows by the Case 1 that there exists a number $l>0$ such that the
congruence $(n^2)^{x}r^2\equiv s^2 \pmod{((n^2)^l-1)}$ has no solution. Then clearly that under
$t=2l$ the congruence $n^{x}r\equiv s \pmod{u_t}$ has no solution too. So, since $r\neq s$ it
remains to consider the case $s=-r$.

Let us suppose, arguing by contradiction, that for every number $t>0$ the congruence
$n^{x}r\equiv -r \pmod{u_{t}}$ is solvable. By the Case 1 there exists a number $t_{0}$ such
that for any number $t>t_{0}$ the congruence $(n^2)^{x}r\equiv -r\pmod{((n^2)^t-1})$  has no
solution. Therefore, if the number $p$  satisfies the  inequality $2^{p-1}>t_0$, then the
solution $x_0$ of congruence $n^{x}r\equiv -r\pmod{((n^{2^p}-1})$ must be an odd number.

Since the numbers $x_0$ ¨ $2^p$ are coprime the greatest common divisor of numbers $n^{x_0}+1$
and $n^{2^p}-1$ is $-n-1$. Consequently, the number $r$ must be divided by any number of form
$$
(-n)^{2^p-1}+(-n)^{2^p-2}+\cdots +(-n) + 1,
$$
where $p>\log_2t_0+1$. But this is impossible since $r\neq 0$. The proof of Proposition 4 is
complete.
\enddemo

\demo{\indent The proof of Theorem 9 in the case $|n|=m$.}

The following statement was actually proved by M.~I.~Kargapolov
[7] but was not stated explicitly:

\proclaim{\indent Proposition 5} Let $C$ be an infinite cyclic
normal subgroup of group $G$. If for every integer $r>0$ the
quotient group $G/C^{r}$ is conjugacy $\Cal F$@-separable then
group $G$ is  conjugacy $\Cal F$@-separable too.\endproclaim

In order to derive from this proposition the conjugacy $\Cal
F$@-separability of groups $G(m,n)$ under $|n|=m$   it is enough
to note that in this case the cyclic subgroup $C=B^{m}$ of group
$G(m,n)$ is infinite and normal in $G(m,n)$. It is clear also that
for any integer $r>0$ the quotient group $$G(m,n)/C^{r}=\langle
a,\,b;\ a^{-1}b^{m}a=b^{\pm m},\,b^{mr}=1\rangle$$ is an
$HNN$@-extension of finite cyclic group. It remains to recall that
by [3] any $HNN$@-extension with finite base group is a conjugacy
$\Cal F$@-separable group.\smallskip

For the completeness of account let me give an outline of proof of
Proposition 5.

So, let $G$ be a group with infinite cyclic normal subgroup $C$
(generated by ele\-ment~$c$) such that for every integer $r>0$ the
quotient group $G/C^{r}$ is conjugacy $\Cal F$@-separable. To
prove that group $G$ is conjugacy $\Cal F$@-separable it is enough
to show that for any elements $f$ and $g$ of group $G$ which are
not conjugate in $G$ there exists an integer $r>0$ such that
elements $f$ and $g$ are not conjugate modulo subgroup~$C^{r}$.

Since in the case when elements $f$ and $g$ are not conjugate
modulo subgroup~$C$ we can put $r=1$, it remains to consider the
case when for some integer $k$ element $f$ is conjugate with
element $gc^{k}$. Obviously, it is sufficient to prove that for
some integer $r>0$ elements $g$ and $gc^{k}$ are not conjugate
modulo subgroup~$C^{r}$.  In order to make this let us introduce
the set of integers $$U=\bigl\{n\in\Bbb Z\bigm |(\exists x\in
G)(x^{-1}gx=gc^{n})\bigr\}$$ and its subset $$V=\bigl\{n\in\Bbb
Z\bigm |(\exists x\in G)(x^{-1}gx=gc^{n}\wedge xc=cx)\bigr\}.$$ It
is easy to see that $V$ is a subgroup of additive group $\Bbb Z$
of integers and if $U\neq V$ then $U$ is the union of $\Bbb Z$ and
some another coset $\Bbb Z+n_{0}$. Note that since elements $g$
and $gc^{k}$ are not conjugate in $G$ the integer $k$ does not
belong to $U$.

Now, for some integer $m\geqslant 0$ we must have $V=m\Bbb Z$. It
is asserting that if $m>0$ then we can put $r=m$, i.~e. elements
$g$ and $gc^{k}$ are not conjugate modulo subgroup~$C^{m}$.
Indeed, if, on the contrary, for some element $x\in G$ and for
some integer $s$ we have $x^{-1}gx=gc^{k+ms}$, then the integer
$k+ms$ belong to $U$ and therefore $k\in U$ but this is
impossible. If $m=0$ then $U=\bigl\{0\bigr\}$ or $U=\bigl\{0,
n_{0}\bigr\}$. If $U=\bigl\{0\bigr\}$ then let $r$ be any positive
integer that does not divide $k$ and if $U=\bigl\{0, n_{0}\bigr\}$
then let $r$ be any positive integer that does not divide both
integers $k$ and $k-n_{0}$. It is clear that then for any integer
$s$ the integer $k+rs$ does not belong to $U$, i.~e. elements $g$
and $gc^{k}$ are not conjugate modulo subgroup~$C^{r}$.

The proof of Theorem 9 is complete.\enddemo

In connection with Theorem 9, the question naturally arises, if
$\pi$ is a set of primes, will the group $G(m,n)$, which is $\Cal
F_{\pi}$@-residual,  be conjugacy $\Cal F_{\pi}$@-residual?
 Above results (Theorem 2 and
Corollary from Theorem 5) exhibit the existence of 1- and
2-elements sets $\pi$ of prime numbers such that the group
$G(1,n)$ is $\Cal F_{\pi}$@-residual. Nevertheless, for the
property to be conjugacy $\Cal F_{\pi}$@-separable is valid the

\proclaim{\indent Theorem 10} {\rm (see [5])} If $n\neq\pm 1$ then
for any set $\pi$ consisting of two prime numbers the group
$G(1,n)$ is not conjugacy $\Cal F_{\pi}$@-separable.\endproclaim

Thus, for any integer $n\neq\pm 1$ there exists a set $\pi$ of
prime numbers such that the group $G(1,n)$ is $\Cal
F_{\pi}$@-residual but is not conjugacy $\Cal F_{\pi}$@-separable.
By contrast, when $|n|=m$, we have:

\proclaim{\indent Theorem 11} {\rm (see [20])} For any set $\pi$
of prime numbers and for any group $G(m,n)$, where $|n|=m$, if
group $G(m,n)$ is $\Cal F_{\pi}$@-residual then it is conjugacy
$\Cal F_{\pi}$@-sepa\-rable.
\endproclaim
\bigskip

\centerline{\bf 3. Subgroup separability of BS-groups}
\medskip

It is well known and easily to see that if $|n|>1$ then in group
$G(1,n)$ the cyclic subgroup $B$ generated by element  $b$ is not
$\Cal F$@-separable. Indeed, element $g=aba^{-1}$ does not belong
to $B$ since in $HNN$-extension $G(1,n)$ it is reduced of length
2. Let $N$ be a finite index normal subgroup of group $G(1,n)$ and
let $r$ be the order of element $b$ modulo $N$. Since elements $b$
and $b^{n}$ are conjugate and therefore have the same order modulo
$N$, the integers $r$ and $n$ are coprime. Hence there exists an
integer $k$ such that $nk\equiv 1\pmod r$ and therefore,
$g=aba^{-1}\equiv ab^{nk}a^{-1}=b^{k}\pmod N$. Thus, element $g$
belongs to subgroup $BN$ for every normal subgroup $N$ of finite
index of group $G(1,n)$ and hence subgroup $B$ is not $\Cal
F$@-separable. Remark that, on the other hand, an arbitrary
non-cyclic finitely generated subgroup of group $G(1,n)$ is of
finite index  and therefore is $\Cal F$@-separable.

In the case $|n|=m$ the situation again appears to be more
definite:

\proclaim{\indent Theorem 12} If $|n|=m$ then the group $G(m,n)$
is subgroup separable. \endproclaim

It should be noted that in the case when $n=m$ this assertion was
long known by the result of [15], which states that any
one-relator group with non-trivial center is subgroup separable.
In general this Theorem was recently proved in [16].

\bigskip

\centerline{\bf References}
\medskip

\item{1.}
{\it Azarov D.~N.} On the virtual residuality of Baumslag --
Solitar groups by finite p@-groups // Modelirovanie i Analyz
Informatsionnykh Sistem. 2013. Vol.~20, No~1. P.~116--123

\item{2.}
{\it Baumslag G., Solitar D.} Some two-generator one-relator
non-Hop\-fian groups // Bull. Amer. Math. Soc. 1962. Vol.~68.
P.~199--201.

\item{3.}
{\it Dyer J.} Separating conjugates in amalgamating free products
and HNN-extensions // J. Austral. Math. Soc. 1980. Vol.~29. No.~1.
P.~35--51.

\item{4.}
{\it Gruenberg K.~W.} Residual properties of infinite soluble
groups // Proc. London Math. Soc. (3) 1957. Vol.~7. P.~29--62.

\item{5.}
{\it Ivanova E.~A., Moldavanskii D.~I.} On the conjugacy
separability of solvable Baumslag -- Solitar groups // Vestnik
Ivanovo State University. Ser. \lq\lq Natural, Social Sci.\rq\rq
2011. Issue~2. P.~129--136 (Russian).

\item{6.} {\it Ivanova O.~A., Moldavanskii D.~I.} Residuality by a finite
$\pi$-groups of some one-relator groups // Sci. Proc. Ivanovo
State University. Mathematics. 2008. Issue.~6. P.~51--58
(Russian).

\item{7.}
{\it Kargapolov M.~I.} Conjugacy separability of supersolvable
groups // Algebra and Logic. 1967. Vol.~6, No.~1. P.~63--68
(Russian).

\item{8.}
{\it Meskin S.} Nonresidually finite one-relator groups // Trans.
Amer. Math. Soc. 1972. Vol.~164. P.~105--114.

\item{9.}
 {\it Moldavanskii D.~I.} On the residuality of Baumslag -- Solitar
 groups //  Chebyshevskii Sb. 2012. Vol.~13, Issue~1. Publisher of
 Tula State Pedagogical University. P.~110 --114 (Russian).

\item{10.}
{\it Moldavanskii D.~I.} The Intersection of the subgroups of
finite index in Baumslag -- Solitar groups // Math. Notes. 2010.
Vol.~87, No.~1, P.~79--86.

\item{11.}
{\it Moldavanskii D.~I.} The Intersection of the subgroups of
finite $p$-index in Baumslag -- Solitar groups // Vestnik Ivanovo
State University. Ser. \lq\lq Natural, Social Sci.\rq\rq\ 2010.
Issue~2. P.~106--111 (Russian).

\item{12.}
{\it Moldavanskii D.~I.} The isomorphism of Baumslag -- Solitar
groups // Ukrainian Math. J. 1991. Vol.~43. No~12. P.~1684--1686
(Russian).

\item{13.}
{\it Moldavanskii D.~I.} Residuality by a finite $p$-groups of
$HNN$-extensions // Vestnik Ivanovo State University.  2000.
Issue~3. P.~129--140 (Russian).

\item{14.}
{\it Moldavanskii D.~I., Kravchenko L.~V., Frolova E.~N.}
// Conjugacy separability of some one-relator groups
// Algorithmic Problems in Group and Semigroup Theory. Edition of
Tula pedagogical inst. 1986. P.~81--91 (Russian).

\item{15.}
{\it Moldavanskii D.~I., Timofeeva L.~V.} Finitely generated
subgroups of one-relator group with non-trivial center are
finitely separable // Izvestiya Vuzov. Mathematics. 1987.
Issue~12. P.~58--59 (Russian).

\item{16.}
{\it Moldavanskii D., Uskova A.} On the finitely separability of
subgroups of generalized free products // arXiv: 1308.3955.
[math.GR].

\item{17.}
{\it Neumann B.~H.} An assay on free products of groups with
amalga\-mations // Phil\. Trans\. Royal Soc\. of London. 1954.
Vol.~246.  P.~503--554.

\item{18.}
{\it Senkevich O.~E.} Conjugacy separability of some
HNN-extensions of groups // Vestnik Ivanovo State University.
2006. Issue~3. P.~133--146 (Russian).

\item{19.} {\it Sokolov E.~V.} Conjugacy separability of descending HNN-extensions
of finitely generated Abelian groups. Math. Notes. 2005. Vol.~78,
No.~5, P.~696--708.

\item{20.}
{\it Varlamova I.~A., Moldavanskii D.~I.} On the residual
finiteness of Baumslag -- Solitar groups // Vestnik Ivanovo State
University. Ser. \lq\lq Natural, Social Sci.\rq\rq 2012. Issue~2.
P.~107--114 (Russian).

\item{21.}
{\it Wong P.~C., Tang C.~K.} Conjugacy separability of certain
$HNN$ extensions // Algebra Colloquium 5:1. 1998.
P.~25--31.\bigskip

Ivanovo State University \smallskip

{\it E-mail address}: moldav\@mail.ru

\end